
\input amstex
\documentstyle{amsppt}
\loadbold
\loadmsbm


\def\bk{{\boldkey K}} 
\def\be{{\boldkey E}} 

\def\ns{\operatorname{ns}}  

\def\Im{\hbox{\rm Im}}

\overfullrule=0pt
\magnification =1100
\vsize 9.25 truein
\hsize 6.5 truein
\voffset -.30in
\TagsOnRight

\document

\topmatter
\title Hankel determinants of Eisenstein series
\endtitle 
\author  Stephen C. Milne \endauthor
\thanks  S. C. Milne was partially supported by 
National Security Agency grants  
MDA 904-97-1-0019 and MDA904--99--1--0003
\endthanks
\affil The Ohio State University \endaffil
\address Department of Mathematics, The Ohio State University, 
Columbus, Ohio, 43210
\endaddress
\email milne\@math.ohio-state.edu 
\endemail
\date 13 September 2000 (revised 2 November 2000)\enddate
\leftheadtext{STEPHEN C. MILNE}
\rightheadtext{Hankel determinants of Eisenstein series}
\abstract
In this paper we prove Garvan's conjectured formula for the square of the 
modular discriminant $\Delta$ as a $3$ by $3$ Hankel determinant
of classical Eisenstein series $E_{2n}$.  We then obtain similar
formulas involving minors of Hankel determinants for $E_{2r}\Delta^m$, for
$m=1,2,3$ and $r=2,3,4,5,7$, and $E_{14}\Delta^4$. We next use Mathematica to
discover, and then the standard structure theory of the ring of modular forms, to
derive the general form of our infinite family of formulas extending the classical
formula for $\Delta$ and Garvan's formula for $\Delta^2$.  This general formula
expresses the $n\times n$ Hankel determinant $\det(E_{2(i+j)}(q))_{1\leq i,j\leq
n}$ as the product of $\Delta^{n-1}(\tau)$, a homogeneous polynomial in $E_4^3$ and
$E_6^2$, and if needed, $E_4$.  We also include a simple verification proof of
the classical $2$ by $2$ Hankel determinant formula for $\Delta$.  This proof
depends upon polynomial properties of elliptic function parameters from Jacobi's
Fundamenta Nova. The modular forms approach provides a convenient explanation
for the determinant identities in this paper. 
\endabstract 
\endtopmatter

\head 1. Introduction\endhead

In this paper we prove Garvan's conjectured formula \cite{11} for the square of 
the modular discriminant $\Delta$ as a $3$ by $3$ Hankel determinant
of classical Eisenstein series $E_{2n}$.  We then obtain similar
formulas involving minors of Hankel determinants for $E_{2r}\Delta^m$, for 
$m=1,2,3$ and $r=2,3,4,5,7$, and $E_{14}\Delta^4$. We next use Mathematica
\cite{37} to discover, and then the modular forms approach of \cite{31, pp.
88--93}, as outlined in \cite{5}, to derive the general form of our infinite family
of formulas extending the classical formula for $\Delta$ and Garvan's formula for 
$\Delta^2$.  This general formula expresses the $n\times n$ Hankel determinant
$\det(E_{2(i+j)}(q))_{1\leq i,j\leq n}$ as the product of
$\Delta^{n-1}(\tau)$, a homogeneous polynomial in $E_4^3$ and
$E_6^2$, and if needed, $E_4$.  We also include a simple verification proof of the
classical formula for $\Delta$ in (1.5) below.  This proof depends upon polynomial
properties of elliptic function parameters from Jacobi's Fundamenta Nova \cite{16}.
The modular forms approach provides a convenient explanation
for the determinant identities in this paper.

The modular discriminant $\Delta$ is defined in \cite{2, Entry 12, pp. 326}
and \cite{28, Eqn. (6.1.11), pp. 196} by means of the following definition.
\definition{Definition 1.1} Let $q:=\exp(2\pi i\tau)$, 
where $\tau$ is in the upper half-plane $\Cal H$.  We then have 
$$\Delta(\tau)\equiv \Delta(q):=q\prod_{r=1}^\infty (1-q^r)^{24}.\tag 1.1$$
\enddefinition  
 
The Fourier expansions of the classical Eisenstein series
$E_n(\tau)$ as given by \cite{2, pp. 318} and
\cite{28, pp. 194--195} are determined by the following definition.
\definition{Definition 1.2} Let $q:=\exp(2\pi i\tau)$, 
where $\tau$ is in the upper half-plane $\Cal H$, and take 
$y:=\Im(\tau)>0$.  Let $n=1,2,3,\cdots$.  We then have 
$$\spreadlines{6 pt}\allowdisplaybreaks\align 
E_2(\tau)\equiv E_2(q):=&1-24
\sum\limits_{r=1}^{\infty}
{rq^r\over {1-q^r}} - {3\over {\pi y}},\tag 1.2\cr
\kern -13.35 em \text{and for } n\geq 2,
\kern 13.35 em\quad&\cr
E_{2n}(\tau)\equiv E_{2n}(q):=&1-{\frac{4n}{B_{2n}}}
\sum\limits_{r=1}^{\infty}
{r^{2n-1}q^r\over {1-q^r}},\tag 1.3\cr
\endalign$$
with the $B_{2n}$ the Bernoulli numbers defined in \cite{7, pp. 48--49} 
by 
$${t\over {e^t-1}}:=\sum\limits_{n=0}^{\infty}
B_n{t^n\over n!}, 
\qquad \text{ for}\  |t|<2\pi.\tag 1.4$$
\enddefinition  

The fundamental classical formula for the modular discriminant
$\Delta$ is provided by the following theorem.
\proclaim{Theorem 1.3} Let $q:=\exp(2\pi i\tau)$, 
where $\tau$ is in the upper half-plane $\Cal H$.  Let $\Delta(\tau)$ 
and $E_{2n}\equiv E_{2n}(q)$ be determined by \hbox{\rm Definitions 1.1} and 
\hbox{\rm 1.2},  respectively. Then, for $|q|<1$, 
$$\Delta(\tau) = {\tfrac{1}{1728}}(E_4^3 - E_6^2).\tag 1.5$$
\endproclaim

Early elliptic function references for (1.5) are \cite{14}, \cite{15, pp. 561}, 
\cite{18, Eqns. (1) and (2), pp. 154}, \cite{25}, and \cite{26, pp. 27}.  (Both 
Hurwitz and Molin replace $q$ by $q^2$.) All of these authors contributions refer
to earlier background developments in \cite{8, 9, 29}.  The chapter notes in
\cite{6, pp. 95, pp. 136--137} provide an excellent summary of Dedekind's
fundamental \cite{8, Eqn. (3), pp. 281; Eqn. (13), pp. 283; Eqn. (24), pp. 285},
Dedekind's comments in \cite{9} on Riemann's work in \cite{29}, and the remarks of
Fricke \cite{10} and Molin \cite{26, pp. 28} on related methods of Jacobi and
Hermite.  The classical elliptic function methods for proving (1.5) are discussed
in \cite{6, pp. 58--72}, \cite{13, pp. 409, pp. 481}, and \cite{22, pp. 125--140,
pp. 177}. An elementary proof of (1.5) was later found by Ramanujan in \cite{27,
Sections 5, 7 and 10}.  This and additional work of Ramanujan involving (1.5) is
surveyed in \cite{3, pp. 114--140}, \cite{4, pp. 43--50}, and \cite{34, pp.
1--18}. Two additional elementary proofs of (1.5) are described in \cite{31, pp.
95--96}.  Recent references for (1.5) are \cite{2, Entry 12(i), pp. 326}, \cite{6,
Theorem 7, pp. 71}, \cite{28, Eqn. (6.1.14), pp. 197}, \cite{30, Eqn. (8), pp. 55;
Theorem 8, pp. 70}, and \cite{31, Eqn (42), pp. 95}.  (The work \cite{19}
is a very useful introduction to \cite{30}.)  Additional applications of (1.5) also
appear in \cite{1, 28, 31}. 

After seeing \cite{23, Theorems 2.1 and 2.2, pp. 15006} and an early version of
\cite{24, Theorems 1.5, 1.6, 5.3--5.6} Garvan \cite{11} observed via
the well-known relation $E_4^2=E_8$ that (1.5) immediately  becomes the Hankel
determinant formula
$$\Delta(\tau) = {\tfrac{1}{1728}}(E_4E_8 - E_6^2) = 
{\tfrac{1}{1728}}\det\vmatrix  E_4 & E_6  \\
  &    \\
E_6  & E_8   
\endvmatrix.\tag 1.6$$ 
He then conjectured the following theorem.
\proclaim{Theorem 1.4~(Garvan)} Let $q:=\exp(2\pi i\tau)$, 
where $\tau$ is in the upper half-plane $\Cal H$.  Let $\Delta(\tau)$ 
and $E_{2n}\equiv E_{2n}(q)$ be determined by \hbox{\rm Definitions 1.1} and 
\hbox{\rm 1.2},  respectively. Then, for $|q|<1$, 
$$\Delta^2(\tau) = -{\tfrac{691}{(1728)^2\cdot 250}}
\det\vmatrix  E_4 & E_6  & E_8  \\
  &    \\
E_6  & E_8   & E_{10}  \\
  &    \\
E_8  & E_{10}   & E_{12}  
\endvmatrix.\tag 1.7$$
\endproclaim
\demo{Proof} Substitute the following three well-known relations from 
\cite{27, Table I., pp. 141}, \cite{28, pp. 195} into the $3$ by $3$ Hankel 
determinant in (1.7).
$$E_8 = E_4^2,\qquad E_{10} = E_4E_6,\qquad 
E_{12} = \tfrac{441}{691}E_4^3 + \tfrac{250}{691}E_6^2.\tag 1.8$$
Simplifying, factoring, and applying (1.5) then gives the $\Delta^2(\tau)$ on the
left-hand-side of (1.7).
\qed\enddemo

For $\Delta^n(\tau)$ with $n>2$, formulas analogous to (1.6) and (1.7) 
generally require a suitable $n+1$ by $n+1$ determinant on the right-hand-side
and an additional polynomial factor in $E_4^3$ and $E_6^2$ on the left-hand-side.  
This extra polynomial factor can often be simplified by relations such as (1.8). 

We organize the rest of our paper as follows. In Section 2 we first apply
recursive methods to obtain our $19$ determinental formulas expressing small
powers of $\Delta(\tau)$, multiplied by a single Eisenstein series, as a suitable
constant times a certain minor of a Hankel determinant of the $E_{2r}$.  These
formulas were motivated by Ramanujan's consideration of $E_{2r}\Delta$, for
$r=2,3,4,5,7$, in \cite{27, Section 16}, and the discussion of these
$E_{2r}\Delta$ in \cite{33, pp. 302}.  The minors here were initially motivated by
the $n$ by $n$ minors of the $n+1$ by $n+1$ Hankel determinants as discussed in
\cite{17, pp. 244--250}.  We next use Mathematica \cite{37} to discover, and
then the modular forms approach of \cite{31, pp. 88--93}, as outlined in
\cite{5}, to derive the general form of our infinite family
of formulas extending (1.6) and (1.7) that involve $\Delta^{n-1}(\tau)$ and an  
$n$ by $n$ Hankel determinant of the $E_{2r}$.

In Section 3 we follow Jacobi's analysis in \cite{16, Section 42}
and utilize the Fourier series for the Jacobi elliptic function $\ns^2$ to write
down a formula for the Eisenstein series $E_{2n}$, for $n\geq 2$. We then apply
\cite{16, Eqn. (2.), Section 36} to put together a simple verification proof of
the classical formula for $\Delta$ in (1.5).  

Symmetry properties of the coefficients in the Maclaurin series expansion of
$\ns^2$ strongly suggest that formulas such as (3.8) and (3.9) in Theorem 3.1 will
be useful in a further study of the determinental formulas in Section 2.  

\head 2. Additional determinental formulas involving powers of $\Delta$
\endhead

Our $19$ determinental formulas in Theorem 2.3 involving small powers of
$\Delta(\tau)$, and the infinite families of identities in Theorem 2.5 are partly
motivated by the determinants in the following definition. 
\definition{Definition 2.1}
Let $\left\{c_{\nu}\right\}_{\nu=1}^\infty$ be a
sequence in ${\Bbb C}^\times$, and let
$m,n=1,2,3,\cdots$.  We take $H_n^{(1)}$ and 
$\chi_n^{(m)}$ to be the determinants of 
$n\times n$ square matrices  
$$\spreadlines{8 pt}\allowdisplaybreaks\align
	H_n^{(1)}\equiv H_n^{(1)} (\{c_\nu\}) 
&:= \det \pmatrix
		c_1 & c_{2} & \ldots & c_{n-1} & c_{n} \\
		c_{2} & c_{3} & \ldots & c_{n} &c_{n+1} \\
		\vdots & \vdots & \ddots & \vdots & \vdots \\
		c_{n} & c_{n+1} & \ldots & c_{2n-2} & c_{2n-1}
	\endpmatrix,\tag 2.1\cr
	\chi_n^{(m)}\equiv \chi_n^{(m)}(\{c_\nu\}) 
&:= \det \pmatrix
		c_1 & c_2 & \ldots & c_{n-m} & c_{n-m+2} & \ldots & c_{n+1} \\
		c_2 & c_3 & \ldots & c_{n-m+1} & c_{n-m+3} & \ldots & c_{n+2} \\
		\vdots & \vdots & \ddots & \vdots & \vdots & \ddots & \vdots \\
		c_n & c_{n+1} & \ldots & c_{2n-m-1} & c_{2n-m+1} & \ldots & c_{2n}
	\endpmatrix.\tag 2.2\cr
\endalign$$
\noindent
The matrix for $\chi_n^{(m)}$ is obtained from the 
matrix for $H_{n+1}^{(1)}$ by deleting the $(n-m+1)$-st 
column and the last row.  Others denote 
$\chi_n^{(1)}$ by $\chi_n$.  We also have $H_0^{(1)} = 1$ and 
$\chi_n^{(m)} = 0$, if $n<m$.  Formally, $\chi_n^{(0)} = H_{n}^{(1)}$.
\enddefinition  

Applications of the Hankel determinants $H_{n}^{(1)}$ and determinants 
$\chi_n^{(m)}$ to continued fractions and orthogonal polynomials are
discussed in \cite{17, pp. 244--250}.  An excellent survey of the literature
on Hankel determinants can be found in Krattenthaler's 
summary in \cite{20, pp. 20--23 ; pp. 46--48}. 

Our eventual aim is to study the combinatorial and geometrical implications of 
the determinants $H_n^{(1)}(\{E_{2(\nu +1)}(q)\})$, 
$\chi_n^{(m)}(\{E_{2(\nu +1)}(q)\})$, and additional minors,  
where $E_{2(\nu +1)}(q)$, with $\nu\geq 1$, are the 
Eisenstein series in Definition 1.2.  A starting point is Theorems 2.3 and 2.5
below. 

Motivated by our proof of Theorem 1.4 we first consider Ramanujan's 
\cite{2, Entry 14, pp. 332} recursion for the $E_{2r}$ given by the following
theorem.
\proclaim{Theorem 2.2~(Ramanujan)} Let $q:=\exp(2\pi i\tau)$, 
where $\tau$ is in the upper half-plane $\Cal H$. Let $E_{2n}(\tau)$ be  
determined by \hbox{\rm Definition 1.2}. For convenience, with $n>1$, define 
$S_{2n}$ by 
$$S_{2n}:= (-1)^{n-1}{\frac{B_{2n}}{4n}}E_{2n}(\tau).\tag 2.3$$
If $n$ is an even integer exceeding $4$ then 
$$\spreadlines{8 pt}\allowdisplaybreaks\align
-{\tfrac{(n+2)(n+3)}{2n(n-1)}}S_{n+2} = &-20{n-2\choose 2}S_4S_{n-2}\cr
&+\sum\limits_{r=1}^{[(n-2)/4]}\kern -.9 em {}^{'}\kern .6 em
{n-2\choose 2r}\biggl\{(n+3-5r)(n-8-5r)\cr &\kern 6 em
-5(r-2)(r+3)\biggr\}S_{2r+2}S_{n-2r},\tag 2.4\cr
\endalign$$
where the prime on the summation sign indicates that if $(n-2)/4$ is an integer,
then the last term of the sum is to be multiplied by $\frac{1}{2}$.
\endproclaim

Substituting $n= 6, 8, 10$ into (2.4) yields the relations in (1.8), and 
setting $n=12$ in (2.4) leads to the well-known relation 
$$E_{14} = E_4^2E_6.\tag 2.5$$

We next utilize (1.5), Theorem 2.2, and Mathematica \cite{37} to derive our
determinental formulas for $E_{2r}\Delta^m$, for $m=1,2,3$ and $r=2,3,4,5,7$, and 
$E_{14}\Delta^4$. In each of the $19$ identities below we first used (2.4) to
write all the $E_{2r}$ in the  determinants on the right-hand-sides as polynomials
in $E_4$ and $E_6$. Simplifying, factoring, applying (1.5), and then referring to
(1.8) and (2.5) as needed yielded the left-hand-side of  each identity.  
We have the following theorem.
\proclaim{Theorem 2.3} Let $q:=\exp(2\pi i\tau)$, 
where $\tau$ is in the upper half-plane $\Cal H$.  Let $\Delta(\tau)$ 
and $E_{2n}\equiv E_{2n}(q)$ be determined by \hbox{\rm Definitions 1.1} and 
\hbox{\rm 1.2},  respectively. Then, for $|q|<1$, 
$$\spreadlines{10 pt}\allowdisplaybreaks\align 
E_4\Delta(\tau) = \ & -{\tfrac{691}{1728\cdot 250}}
\det\vmatrix  E_4 & E_8  \\
  &    \\
E_8  & E_{12}    
\endvmatrix,\tag 2.6\cr
E_6\Delta(\tau) = \ & {\tfrac{691}{1728\cdot 250}}
\det\vmatrix  E_4 & E_{12}  \\
  &    \\
E_6  & E_{14}    
\endvmatrix,\tag 2.7\cr
E_8\Delta(\tau) = \ & {\tfrac{3617}{1728\cdot 3\cdot 7^2\cdot 11}}
\det\vmatrix  E_4 & E_{10}  \\
  &    \\
E_{10}  & E_{16}    
\endvmatrix\tag 2.8\cr
 = \ & {\tfrac{691}{1728\cdot 441}}
\det\vmatrix  E_8 & E_{10}  \\
  &    \\
E_{10}  & E_{12}    
\endvmatrix,\tag 2.9\cr
E_{10}\Delta(\tau) = \ & {\tfrac{691\cdot 43867}
{1728\cdot 2\cdot 3^2\cdot 5^4\cdot 7^2\cdot 13}}
\det\vmatrix  E_4 & E_{12}  \\
  &    \\
E_{10}  & E_{18}    
\endvmatrix\tag 2.10\cr
 = \ & {\tfrac{691}{1728\cdot 250}}
\det\vmatrix  E_8 & E_{12}  \\
  &    \\
E_{10}  & E_{14}    
\endvmatrix,\tag 2.11\cr
E_{14}\Delta(\tau) = \ & {\tfrac{691\cdot 593\cdot 131}
{1728\cdot 2\cdot 3\cdot 5^3\cdot 7^2\cdot 11\cdot 13}}
\det\vmatrix  E_4 & E_{14}  \\
  &    \\
E_{12}  & E_{22}    
\endvmatrix\tag 2.12\cr
 = \ & -{\tfrac{691\cdot 3617}
{1728\cdot 2\cdot 3\cdot 5^3\cdot 7^2\cdot 13}}
\det\vmatrix  E_{10} & E_{14}  \\
  &    \\
E_{12}  & E_{16}    
\endvmatrix,\tag 2.13\cr
E_4\Delta^2(\tau) = \ & {\tfrac{691^2}{(1728)^2\cdot (21)^2\cdot 250}}
\det\vmatrix  E_4 & E_8  & E_{10}  \\
  &    \\
E_6  & E_{10}   & E_{12}  \\
  &    \\
E_8  & E_{12}   & E_{14}  
\endvmatrix,\tag 2.14\cr
E_6\Delta^2(\tau) = \ & -{\tfrac{691^2}{(1728)^2\cdot (250)^2}}
\det\vmatrix  E_6 & E_8  & E_{10}  \\
  &    \\
E_8  & E_{10}   & E_{12}  \\
  &    \\
E_{10}  & E_{12}   & E_{14}  
\endvmatrix,\tag 2.15\cr
E_8\Delta^2(\tau) = \ & -{\tfrac{(691)^2\cdot 3617}
{(1728)^2\cdot 2^3\cdot 3\cdot 5^3\cdot 7^2\cdot 467}}
\det\vmatrix  E_6 & E_8  & E_{12}  \\
  &    \\
E_8  & E_{10}   & E_{14}  \\
  &    \\
E_{10}  & E_{12}   & E_{16}  
\endvmatrix,\tag 2.16\cr
E_{10}\Delta^2(\tau) = \ & {\tfrac{(691)^2\cdot 3617}
{(1728)^2\cdot 2^2\cdot 3\cdot 5^6\cdot 7^2\cdot 13}}
\det\vmatrix  E_6 & E_{10}  & E_{12}  \\
  &    \\
E_8  & E_{12}   & E_{14}  \\
  &    \\
E_{10}  & E_{14}   & E_{16}  
\endvmatrix,\tag 2.17\cr
E_{14}\Delta^2(\tau) = \ & -{\tfrac{(691)^2\cdot 3617\cdot 43867}
{(1728)^2\cdot 2^6\cdot 3\cdot 5^3\cdot 7^2\cdot 97\cdot 7213}}
\det\vmatrix  E_8 & E_{10}  & E_{14}  \\
  &    \\
E_{10}  & E_{12}   & E_{16}  \\
  &    \\
E_{12}  & E_{14}   & E_{18}  
\endvmatrix,\tag 2.18\cr
E_4\Delta^3(\tau) = \ & -{\tfrac{(691)^3\cdot 3617}
{(1728)^3\cdot 2^4\cdot 3\cdot 5^6\cdot 7^2\cdot 467}}
\det\vmatrix  E_4 & E_6 & E_8  & E_{10}  \\
  &    \\
E_6 & E_8  & E_{10}   & E_{12}  \\
  &    \\
E_8 & E_{10}  & E_{12}   & E_{14}  \\
  &    \\
E_{10} & E_{12}  & E_{14}   & E_{16}  
\endvmatrix,\tag 2.19\cr
E_6\Delta^3(\tau) = \ & -{\tfrac{(691)^3\cdot 43867}
{(1728)^3\cdot 2^5\cdot 3^2\cdot 5^6\cdot 7^2\cdot 131}}
\det\vmatrix  E_4 & E_6 & E_8  & E_{12}  \\
  &    \\
E_6 & E_8  & E_{10}   & E_{14}  \\
  &    \\
E_8 & E_{10}  & E_{12}   & E_{16}  \\
  &    \\
E_{10} & E_{12}  & E_{14}   & E_{18}  
\endvmatrix,\tag 2.20\cr
E_8\Delta^3(\tau) = \ & {\tfrac{(691)^3\cdot (3617)^2}
{(1728)^3\cdot 2^4\cdot 3^2\cdot 5^6\cdot 7^4\cdot 13\cdot 467}}
\det\vmatrix  E_4 & E_6 & E_{10}  & E_{12}  \\
  &    \\
E_6 & E_8  & E_{12}   & E_{14}  \\
  &    \\
E_8 & E_{10}  & E_{14}   & E_{16}  \\
  &    \\
E_{10} & E_{12}  & E_{16}   & E_{18}  
\endvmatrix,\tag 2.21\cr
E_{10}\Delta^3(\tau) = \ & -{\tfrac{(691)^3\cdot 3617\cdot 43867}
{(1728)^3\cdot 2^7\cdot 3\cdot 5^6\cdot 7^2\cdot 97\cdot 7213}}
\det\vmatrix  E_4 & E_8 & E_{10}  & E_{12}  \\
  &    \\
E_6 & E_{10}  & E_{12}   & E_{14}  \\
  &    \\
E_8 & E_{12}  & E_{14}   & E_{16}  \\
  &    \\
E_{10} & E_{14}  & E_{16}   & E_{18}  
\endvmatrix,\tag 2.22\cr
E_{14}\Delta^3(\tau) = \ & {\tfrac{(691)^3\cdot (3617)^2\cdot 43867\cdot 283
\cdot 617}
{(1728)^4\cdot 2^3\cdot 3^2\cdot 5^6\cdot 7^2\cdot 31\cdot 3503110621}}
\det\vmatrix  E_6 & E_8 & E_{10}  & E_{14}  \\
  &    \\
E_8 & E_{10}  & E_{12}   & E_{16}  \\
  &    \\
E_{10} & E_{12}  & E_{14}   & E_{18}  \\
  &    \\
E_{12} & E_{14}  & E_{16}   & E_{20}  
\endvmatrix,\tag 2.23\cr
E_{14}\Delta^4(\tau) = \ & {\tfrac{(691)^4\cdot (3617)^2\cdot 43867\cdot 
131\cdot 283\cdot 593\cdot 617}
{(1728)^6\cdot 5^9\cdot 7^5\cdot 11\cdot 13\cdot 67\cdot 257\cdot 43721}}
\det\vmatrix  E_4 & E_6 & E_8  & E_{10} & E_{14} \\
  &    \\
E_6 & E_8  & E_{10}   & E_{12} & E_{16} \\
  &    \\
E_8 & E_{10}  & E_{12}   & E_{14} & E_{18} \\
  &    \\
E_{10} & E_{12}  & E_{14}   & E_{16}  & E_{20}  \\
  &    \\
E_{12} & E_{14}  & E_{16}   & E_{18}  & E_{22} 
\endvmatrix.\tag 2.24\cr
\endalign$$
\endproclaim

The determinants in (1.6), (1.7), (2.14), (2.15), (2.19)--(2.22), (2.24) are of 
the form $H_2^{(1)}$, $H_3^{(1)}$, $\chi_3^{(2)}$, $\chi_3^{(3)}$,
$H_4^{(1)}$, $\chi_4$, $\chi_4^{(2)}$, $\chi_4^{(3)}$, and $\chi_5$, 
respectively, with entries $c_{\nu} = E_{2(\nu +1)}(q)$.  The rest are certain
other minors of $H_n^{(1)}(\{E_{2(\nu +1)}(q)\})$, for suitable $n$. 

The determinant evaluations in (2.6)--(2.24) can also be proven, as pointed out
in \cite{5}, by the methods in the modular forms proof of Theorem 2.5 below. 
Being able to apply the first part of this proof leads to a characterization of
a large class of determinants of Eisenstein series which have
evaluations analogous to those in (2.6)--(2.24) and Theorem 2.5. It turns out we
only have to consider minors of $H_n^{(1)}(\{E_{2(\nu +1)}(q)\})$, for suitable
$n$.  We have the following proposition.  
\proclaim{Proposition 2.4}Let $A$ be any $n\times n$ square matrix whose entries
are Eisenstein series $E_{2r}$, with $r\geq 2$.  Suppose there are no repeated
rows or columns. Recall that the weight of a product 
$E_{2r_1}E_{2r_2}\cdots E_{2r_n}$ of Eisenstein series is 
$2(r_1+r_2+\cdots+r_n)$. Then, each term in the $n\times n$ determinant 
$\det A$ has the same weight if and only if $\det A$ is $\pm 1$ times some 
$n\times n$ minor of a Hankel determinant 
$H_m^{(1)}(\{E_{2(\nu +1)}(q)\})
\equiv\det(E_{2(i+j)}(q))_{1\leq i,j\leq m}$ of Eisenstein series, with 
$m\geq n$.    
\endproclaim
\demo{Proof}First, assume that $\det A$ is some $n\times n$ minor of 
$\det(E_{2(i+j)}(q))_{1\leq i,j\leq m}$, with $m\geq n$.  Let the rows and
columns of $A$ be indexed by $\{i_1,\dots i_n\}$ and $\{j_1,\dots j_n\}$, 
respectively.  The weight of the term corresponding to $\sigma$ in 
$$\det A = \sum_{\sigma\in \Cal S_n}\text{sign}(\sigma)
\prod\limits_{r=1}^{n}E_{2(i_r+j_{\sigma(r)})}(q),\tag 2.25$$ 
is $2((i_1+\cdots i_n)+(j_{\sigma(1)}+\cdots j_{\sigma(n)}))=
2((i_1+\cdots i_n)+(j_1+\cdots j_n))$, which is a constant.

Next, for an even more general argument, suppose that each term in the $n\times n$
determinant 
$$\det A \equiv\det(E_{p_{i,j}}(q))_{1\leq i,j\leq n}= 
\sum_{\sigma\in \Cal S_n}\text{sign}(\sigma)
\prod\limits_{r=1}^{n}E_{p_{r,\sigma(r)}}(q),\tag 2.26$$
has the same weight, where the $n^2$ subscripts $p_{i,j}$ are now arbitrary
reals. In this setting we take the weight of the term in (2.26) 
corresponding to $\sigma$ to be 
$$p_{1,\sigma(1)}+p_{2,\sigma(2)}+\cdots +p_{n,\sigma(n)}.
\tag 2.27$$
We then claim that 
$$p_{i,j}-p_{i,j-1} = p_{1,j}-p_{1,j-1},\quad\text{for}\ i=1,2,\cdots,n
\quad\text{and}\ j=2,3,\cdots,n,\tag 2.28$$
and
$$p_{i,j}-p_{i-1,j} = p_{i,n}-p_{i-1,n},\quad\text{for}\ i=2,3,\cdots,n
\quad\text{and}\ j=1,2,\cdots,n.\tag 2.29$$

To obtain (2.28) and (2.29) first consider the $2\times 2$ submatrix
$$\vmatrix  E_{p_{i-1,j-1}}(q) & E_{p_{i-1,j}}(q)  \\
  &    \\
E_{p_{i,j-1}}(q)  & E_{p_{i,j}}(q)    
\endvmatrix,\tag 2.30$$
where $i,j=2,3,\cdots,n$.  Keeping in mind (2.30), there are at least two terms in
(2.26) of the form $B\cdot E_{p_{i-1,j-1}}(q)E_{p_{i,j}}(q)$ and 
$-B\cdot E_{p_{i,j-1}}(q)E_{p_{i-1,j}}(q)$.  Equating the weights of these 
two terms and simplifying, gives 
$$p_{i-1,j-1}+p_{i,j} = p_{i,j-1}+p_{i-1,j},\quad\text{for}\ i,j=2,3,\cdots,n.
\tag 2.31$$
By rewriting (2.31) in two ways, we have 
$$p_{i,j}-p_{i,j-1} = p_{i-1,j}-p_{i-1,j-1},\quad\text{for}\ i,j=2,3,\cdots,n,
\tag 2.32$$
and
$$p_{i,j}-p_{i-1,j} = p_{i,j-1}-p_{i-1,j-1},\quad\text{for}\ i,j=2,3,\cdots,n.
\tag 2.33$$
Equation (2.28) is immediate from (2.32) by fixing $j$ and varying $i$, while 
(2.29) is immediate from (2.33) by fixing $i$ and varying $j$.  

By permuting the columns and then the rows of the matrix $A$ in (2.26), and
factoring out a $-1$ if necessary, we can assume that the
differences in the right-hand-sides of (2.28) and (2.29) are strictly positive. 

It is now clear that if we take the $n^2$ subscripts $p_{i,j}$ to be even integers 
greater than $2$, then $\det A$ is $\pm 1$ times some 
$n\times n$ minor of a Hankel determinant $\det(E_{2(i+j)}(q))_{1\leq i,j\leq m}$
of Eisenstein series, with $m\geq n$.
\qed\enddemo

The simplest application of Proposition 2.4 involves applying the 
modular forms approach of \cite{31, pp. 88--93}, as outlined in \cite{5}, to 
express the $n\times n$ Hankel determinant $H_n^{(1)}(\{E_{2(\nu +1)}(q)\})$
as the product of $\Delta^{n-1}(\tau)$, a homogeneous polynomial in $E_4^3$ and
$E_6^2$, and if needed, $E_4$.  One of our original motivations for studying the
determinants $H_n^{(1)}(\{E_{2(\nu +1)}(q)\})$ was to note that 
$$E_4 = H_1^{(1)}(\{E_{2(\nu +1)}(q)\}) = \det|E_4|,\tag 2.34$$
and recall the $n = 2, 3, 4$ cases in (1.6), (1.7), and (2.19), respectively.  
Since $E_4$ only appears as a factor in the left-hand-sides of (2.34) and (2.19),
it was natural to split $n$ up into the classes $(\text{mod}\ 3)$ given by 
$n = 3r+1$, $3r+2$, and $3r+3$, with $r=0,1,2,\cdots$.

Mathematica \cite{37} computations of $H_n^{(1)}(\{E_{2(\nu +1)}(q)\})$ up to 
$n = 10$, analogous to those of (1.7) and (2.19), first led to the discovery of 
the general form of the following three infinite families of formulas in 
(2.35)--(2.37) below.  We prove the evaluations in (2.35)--(2.37) by appealing to
the standard structure theory of the ring of modular forms in \cite{31, pp.
88--93}.  The proof here is a detailed rewriting of the original proof supplied by
Borcherds in \cite{5}.  We have the following theorem.  
\proclaim{Theorem 2.5} Let $q:=\exp(2\pi i\tau)$, 
where $\tau$ is in the upper half-plane $\Cal H$.  Let $\Delta(\tau)$ 
and $E_{2m}\equiv E_{2m}(q)$ be determined by \hbox{\rm Definitions 1.1} and 
\hbox{\rm 1.2},  respectively. Then, for $|q|<1$, 
$$E_4\Delta^{3r}(\tau)\cdot P_{3r(r-1)/2}(E_4^3,E_6^2) 
= \ d_r
H_{3r+1}^{(1)}(\{E_{2(\nu +1)}(q)\}),\quad\text{for}\quad r=0,1,2,\cdots, 
\tag 2.35$$

$$\Delta^{3r+1}(\tau)\cdot Q_{r(3r-1)/2}(E_4^3,E_6^2) 
= \ e_r 
H_{3r+2}^{(1)}(\{E_{2(\nu +1)}(q)\}),\quad\text{for}\quad r=0,1,2,\cdots, 
\tag 2.36$$

$$\Delta^{3r+2}(\tau)\cdot R_{r(3r+1)/2}(E_4^3,E_6^2) 
= \ f_r 
H_{3r+3}^{(1)}(\{E_{2(\nu +1)}(q)\}),\quad\text{for}\quad r=0,1,2,\cdots, 
\tag 2.37$$
where $d_r$, $e_r$, and $f_r$ are constants depending on $r$, and 
$P_n(x,y)$, $Q_n(x,y)$, and $R_n(x,y)$ are homogeneous polynomials in $x$ and $y$
of total degree $n$ (as given above), with integer coefficients, whose monomials
are those in $(x-y)^n$.  
\endproclaim
\demo{Proof}Let $H_n^{(1)}(\{E_{2(\nu +1)}(q)\})$ denote any of the Hankel
determinants in (2.35)--(2.37).  Keeping in mind that $E_{2r}(q)$ has weight $2r$,
it is immediate from Proposition 2.4 that each term in the Hankel determinant 
$H_n^{(1)}(\{E_{2(\nu +1)}(q)\})$ is an entire modular form of fixed weight 
$2n(n+1)$, as is the entire determinant.

Each Eisenstein series $E_{2r}(q)$ in (1.3) is written in \cite{31, Eqn.
(34), pp. 92} as a Maclaurin series in $q$ starting with the terms $1+a_1q$, with 
$a_1\neq 0$.  Thus, subtracting the first row from each of the others, factoring
$q$ out of each of the resulting $n-1$ lower rows, and recalling 
$q:=\exp(2\pi i\tau)$, we find that $H_n^{(1)}(\{E_{2(\nu +1)}(q)\})$
vanishes to order $n-1$ at the cusp $\tau=i\infty$.

The Maclaurin series in $q$ for $\Delta(q)$ in (1.1) starts with the term $q$. 
The function $\Delta(q)\equiv\Delta(\tau)$ is also a cusp form of weight $12$
which vanishes at $\tau=i\infty$.  It follows from \cite{31, Theorem 4
(iii), pp. 88} that $\Delta^{n-1}(\tau)$ divides 
$H_n^{(1)}(\{E_{2(\nu +1)}(q)\})$, and that the quotient is a holomorphic modular
form of weight $2n(n+1)-12(n-1) = 2(n-2)(n-3)$.  

By Corollary 2 of \cite{31, pp. 89} the quotient 
$H_n^{(1)}(\{E_{2(\nu +1)}(q)\})/\Delta^{n-1}(\tau)$ is a linear combination of
monomials 
$$E_4^{\alpha}E_6^{\beta},\tag 2.38$$
with fixed weight 
$$4\alpha + 6\beta = 2(n-2)(n-3),\tag 2.39$$
where $\alpha$ and $\beta$ are nonnegative integers.  To obtain the
left-hand-sides of (2.35)--(2.37) we next utilize the weight 
$2(n-2)(n-3)$ $\text{mod}\ 12$ to simplify (2.38).  

First, as in \cite{31, pp. 90}, let $4\alpha + 6\beta = 12m$, with $m$ a
nonnegative integer.  This gives $2\alpha + 3\beta = 6m$, which when divided by
$2$ and $3$, respectively, implies that $\beta/2=b$ and $\alpha/3=a$ are
nonnegative integers.  This gives 
$$E_4^{\alpha}E_6^{\beta}=E_4^{3a}E_6^{2b},\tag 2.40$$   
which in turn leads to (2.36) and (2.37) where $n=3r+2$ and $3r+3$, respectively.  
The total degree of the polynomials $Q$ and $R$ in (2.36) and (2.37) is immediate
from (2.39) by computing 
$$(4\alpha + 6\beta)/12 = a+b =(n-2)(n-3)/6,\tag 2.41$$
for $n=3r+2$ and $3r+3$, respectively.  

Next, to obtain (2.35), let $4\alpha + 6\beta = 12m+4$.  This gives 
$2\alpha + 3\beta = 6m+2$.  We must have $\alpha\geq 1$ or else $\beta$ is not an
integer.  Dividing by $2$ implies that $\beta/2=b$ is a nonnegative integer. 
Setting $\beta = 2b$ and solving for $\alpha$ gives $\alpha = 3(m-b)+1 := 3c+1$,
with $c$ a nonnegative integer.  We now have 
$$E_4^{\alpha}E_6^{\beta}=E_4\cdot E_4^{3c}E_6^{2b},\tag 2.42$$ 
which in turn leads to (2.35) when $n=3r+1$.  In this case we have 
$$(4\alpha + 6\beta -4)/12 = c+b =[(n-2)(n-3)-2]/6\tag 2.43$$
which equals $3r(r-1)/2$.  
\qed\enddemo

It is an interesting open problem to find a concise combinatorial and/or
analytical description of the coefficients in the polynomials 
$P_{3r(r-1)/2}$, $Q_{r(3r-1)/2}$, and $R_{r(3r+1)/2}$.  

The $r=0$ cases of (2.35), (2.36), and (2.37) are (2.34), (1.6), and (1.7),
respectively.  The first nontrivial case of (2.35) is (2.19).  The degree of 
$P_{3r(r-1)/2}$ in (2.35) is $3{r\choose 2}$, while the degrees of 
$Q_{r(3r-1)/2}$ and $R_{r(3r+1)/2}$ in (2.36) and (2.37), respectively, are the
pentagonal numbers $r(3r\mp 1)/2$ in \cite{32, sequence \bf{M1336}}.  See also
\cite{6, pp. 124}. 

Given (2.14), (2.15), (2.20)--(2.22), (2.24), and the modular forms proof of 
Theorem 2.5, we can obtain more complicated infinite families of identities
analogous to (2.35)--(2.37) for $\chi_n^{(m)}(\{E_{2(\nu +1)}(q)\})$.  

Borcherds also observed in \cite{5} that the method of proof of Theorem 2.5 also 
establishes the determinental identities in Theorem 2.3, up to some constant. 
The space of modular forms of the appropriate weight happens to have dimension
$1$, and is thus spanned by an Eisenstein series.  The argument through equation 
(2.39) is the same, with the right-hand-side of (2.39) replaced by the sum of the
subscripts of the Eisenstein series in the diagonal entries of the given 
$n\times n$ matrix, minus $12(n-1)$.  Call this expression $W_1$. We then look at 
$W_1$ $\text{mod}\ 12$, as before. We find that any holomorphic modular form of 
weight $0$,$2$,$4$,$6$,$8$, or $10$ $\text{mod}\ 12$ is equal to $1$, $E_{14}$, 
$E_4$, $E_6$, $E_8$, or $E_{10}$ times a homogeneous polynomial in $E_4^3$ and 
$E_6^2$.  In the case of Theorem 2.3 we have $W_1 = 4$,$6$,$8$,$10$,$14$.  These
weights correspond to the right spaces of dimension $1$ in our list $\text{mod}\
12$, the homogeneous polynomial is a constant, and we are done. We have to be 
careful in the case of $W_1 =14$.  Here, we start with $W_1 = 4\alpha + 6\beta =
12m+2$, for $m\geq 1$. We end up factoring out $E_4^2E_6$, with the 
homogeneous polynomial having total degree $m-1$. The other cases are simpler, 
use $m\geq 0$, and the homogeneous polynomials all have total degree $m$. 

The above modular forms proofs of Theorems 2.3 and 2.5, combined with the
characterization in Proposition 2.4, lead to only a small number of
determinental identities in which the space of modular forms of the appropriate
weight $W_1$ has dimension $1$. In particular, as soon as $W_1\geq 16$, the
dimension is $\geq 2$, and the homogeneous polynomial in $E_4^3$ and $E_6^2$ is no
longer a constant. All such identities analogous to those in Theorems 1.3, 1.4,
and 2.3 are determined by requiring $W_1\leq 14$.  It is not hard to see that this
is not possible if $n\geq 6$.  The remaining finite number of possible cases for
$n\leq 5$ can be checked directly.  The identities in Theorems 1.3, 1.4, and 2.3
cover the basic types that are possible. Moreover, the identities in (1.6) and
(1.7) are the only ones involving just a power of $\Delta$, and the identity in
(2.24) is unique up to transposition symmetry of the $5\times 5$ determinant. 
There are no other analogous $5\times 5$ determinental identities whose
left-hand-side is of the form $E_{2r}\Delta^4(\tau)$. Thus, (2.24) should be very
interesting.   

Keeping in mind \cite{35, Eqn. (52.6), pp. 201} it is natural to consider the
Hankel determinants $H_n^{(1)}(\{E_{2(\nu +1)}(q)\})$ in which entries 
$E_{2(\nu +1)}(q)$ are replaced by $0$ unless $2(\nu +1)$ satisfies any of a fixed
set of congruence conditions.  (The ``unless'' can also be replaced by
``whenever''). That is, when all entries in certain of the counter diagonals are
$0$.  For example, the condition $E_{2(\nu +1)}(q)\mapsto 0$ unless 
$2(\nu +1)\equiv 0\pmod{6}$, leads to interesting determinant evaluations.
Similarly, the condition $E_{2(\nu +1)}(q)\mapsto 0$ whenever 
$2(\nu +1)\equiv 0\pmod{4}$ or $2(\nu +1)\equiv 0\pmod{6}$ leads to reasonable
determinants.  This is just a small sample of many such possibilities.

\head 3. The Jacobi elliptic function $\ns^2$ and Eisenstein series
\endhead

In this section we follow Jacobi's analysis in \cite{16, Section 42}
and utilize the Fourier series for the Jacobi elliptic function $\ns^2$ to write
down a formula for the Eisenstein series $E_{2n}$, for $n\geq 2$. We then apply
\cite{16, Eqn. (2.), Section 36} to put together a simple verification proof of
the classical formula for $\Delta$ in (1.5). 

We require the Jacobi elliptic function parameter
$$z:={}_2F_1\left[\left.\matrix 
\tfrac{1}{2}\ ,\ \tfrac{1}{2}\\ 1\endmatrix\ \right|
\ k^2\right]  = 2\bk(k)/\pi
\equiv 2\bk/\pi,\tag 3.1$$ with 
$$\bk(k)\equiv \bk := \int_0^1 \frac{dt}
{\sqrt{(1-t^2)(1-k^2 t^2)}}
		= \tfrac{\pi}{2}\ 
{}_2F_1\left[\left.\matrix \tfrac{1}{2}\ ,\ 
\tfrac{1}{2}\\
1\endmatrix\ \right|\ k^2\right]\tag 3.2$$   
the complete elliptic integral of the first kind in 
\cite{21, Eqn. (3.1.3), pp. 51}, and $k$ the modulus.  
We also need the complete elliptic integral of the
second kind 
$$\be(k)\equiv \be := \int_0^1 
\sqrt{\frac{1-k^2 t^2}{1-t^2}}\,dt
		= \tfrac{\pi}{2}\ 
{}_2F_1\left[\left.\matrix \tfrac{1}{2}\ ,\ 
-\tfrac{1}{2}\\
1\endmatrix\ \right|\ k^2\right].\tag 3.3$$
Finally, we take 
$$q:=exp(-\pi  {\bk}(\sqrt{1-k^2})/{\bk}(k))\tag 3.4$$

The classical Fourier expansion for $\ns^2$, which first appeared in 
\cite{16, Eqn. (2.), Section 42; Eqn. IV., Section 44}, is now given by 
$$\ns^2(u,k) = 1-\frac{\be}{\bk}
  +\frac{1}{z^2}\csc^2 \tfrac{u}{z} -
\frac{8}{z^2}
\sum\limits_{n=1}^{\infty}
		\frac{n q^{2n}}{1-q^{2n}} \cos \tfrac{2nu}{z}.
	\tag 3.5$$ 
More recent references for (3.5) include 
\cite{12, Eqn. 27, pp. 913} and \cite{36, Ex. 57, pp. 535}.  

The Fourier expansion in (3.5) may be written as a
double sum by first expanding the $\cos\tfrac{2nu}{z}$ as a Maclaurin series,
interchanging summation, and then simplifying. Next, the Laurent series 
expansion for $\csc^2 \tfrac{u}{z}$ is immediate from differentiating the 
Laurent series for $\cot\theta$ in \cite{7, Ex. 36, pp. 88}.  This analysis
yields 
$$\spreadlines{8 pt}\allowdisplaybreaks\align
\ns^2(u,k) = \ &\frac{1}{u^2}+1-\frac{\be}{\bk}+
\sum\limits_{m=1}^{\infty}B_{2m}{\frac{(-1)^{m+1}2^{2m-1}}{m\cdot z^{2m}}}
{\frac{u^{2m-2}}{(2m-2)!}}\cr  
&-4\kern-.2em\sum\limits_{m=1}^{\infty}
   \frac{(-1)^{m-1}2^{2m-1}}{z^{2m+2}}
    \kern-.25em\left[\sum\limits_{r=1}^{\infty}
		\frac{r^{2m-1}q^{2r}}{1-q^{2r}}\right]
  \kern-.25em\frac{u^{2m-2}}{(2m-2)!},\tag 3.6\cr
\endalign$$
with $B_{2m}$ the Bernoulli numbers defined by (1.4).  

The Laurent series expansion for $\ns^2(u,k)$ is 
$$\ns^2(u,k) = \frac{1}{u^2}+\sum_{m=1}^\infty
		 (ns^2)_m(k^2)\frac{u^{2m-2}}{(2m-2)!},\tag 3.7$$
where $(ns^2)_m(k^2)$ are polynomials in $k^2$, with $k$ the modulus.

In what follows we equate the $q$'s in (3.4) and 
Definition 1.2.  That is $2\tau = i {\bk}(\sqrt{1-k^2})/{\bk}(k)$. 

Keeping in mind Definition 1.2, we find that equating coefficients of $u^{2m-2}$,
for $m\geq 2$, in (3.6) and (3.7) leads to a formula for $E_{2m}(q^2)$. 
Furthermore, applying the Gau{\ss} transformation ($q\mapsto \sqrt{q}$,
$k\mapsto \frac{2\sqrt{k}}{1+k}$, $\bk\mapsto (1+k)\bk$, $z\mapsto (1+k)z$) from
\cite{16, Theorem III, Section 37} to the first formula yields a corresponding
formula for $E_{2m}(q)$. We have the following theorem. 
\proclaim{Theorem 3.1} Let $z:= 2\bk(k)/\pi
\equiv 2\bk/\pi$, as in \hbox{\rm(3.1)}, with $k$ the
modulus.  Let the Bernoulli numbers $B_{2m}$ be defined by \hbox{\rm(1.4)}. 
Take $(ns^2)_m(k^2)$ to be the elliptic function polynomials of $k^2$
determined by \hbox{\rm(3.7)}. Let $q$ be as in \hbox{\rm(3.4)}. Take $E_{2m}(q)$
as in \hbox{\rm Definition 1.2}, with 
$2\tau = i {\bk}(\sqrt{1-k^2})/{\bk}(k)$. Let
$m=2,3,4,\cdots$.  We then have
$$\spreadlines{8 pt}\allowdisplaybreaks\align 
E_{2m}(q^2)&\equiv E_{2m}(2\tau) =\ 1-z^2 + 
{(-1)^{m-1}m\cdot z^{2m+2} \over 2^{2m-1}\cdot B_{2m}}
\cdot (ns^2)_m(k^2),\tag 3.8\cr
E_{2m}(q)&\equiv E_{2m}(\tau) =\ 1-(1+k)^2z^2 \cr
&\kern 6 em + {(-1)^{m-1}m\cdot z^{2m+2}(1+k)^{2m+2}
\over 2^{2m-1}\cdot B_{2m}}
\cdot (ns^2)_m(4k/(1+k)^2).\tag 3.9\cr
\endalign$$
\endproclaim

The $m=2$ and $3$ cases of (3.8) are given by 
$$\spreadlines{8 pt}\allowdisplaybreaks\align 
E_{4}(q^2) =\ &z^4(1-k^2+k^4),\tag 3.10\cr
E_{6}(q^2) =\ &z^6(1+k^2)(1-2k^2)(1-\tfrac{1}{2}k^2).\tag 3.11\cr
\endalign$$

Equations (3.10) and (3.11) appear in \cite{3, Entry 13(i),(ii), pp. 126}. 
Equation (3.10) is also recorded in \cite{4, Eqn. (12.21), pp. 49}, and 
just under equation (6.) of \cite{16, Section 42}.  Jacobi's derivation of his $8$
squares formula in Section 42 of \cite{16} only required him to go as far as 
$E_{4}(q^2)$. 

Our verification proof of (1.5) is a consequence of (3.10), (3.11), and 
equation (2.) of \cite{16, Section 36} written in our notation as
$$\Delta^{1/4}(2\tau)\equiv \Delta^{1/4}(q^2) = \tfrac{1}{2^2}z^3k(1-k^2)^{1/2}
.\tag 3.12$$
The fourth power of (3.12) immediately gives 
$$\Delta(2\tau)\equiv \Delta(q^2) = \tfrac{1}{2^8}z^{12}(1-k^2)^2k^4.\tag 3.13$$
Equation (3.13) also appears in \cite{3, Entry 12(iii), pp. 124} and \cite{4,
Theorem 8.3(iii), pp. 29}.

Substituting (3.10) and (3.11) into the right-hand-side of (1.5), with 
$q$ replaced by $q^2$, simplifying, and obtaining the right-hand-side of 
(3.13) now completes the proof of (1.5).  The $z^{12}$ factored out quickly and 
reduced the proof to a computation involving polynomials of low degree in 
$k^2$.  
 
Our simple verification proof of (1.5) does not seem to have been written down
in the literature before. 

\noindent{\bf Acknowledgment:}~We would like to thank R. Borcherds for pointing out
in \cite{5} how the standard structure theory of the ring of modular forms yields 
the above proof of Theorem 2.5, and also how these same techniques may be used to
give an alternate proof of Theorem 2.3. 



\widestnumber\no{9999} 
\Refs

\ref\key 1
\by T. M. Apostol
\paper Modular Functions and Dirichlet Series in Number Theory
\inbook vol. 41 of Graduate Texts in Mathematics
\publ Springer-Verlag\publaddr New York 
\yr 1976
\endref

\ref\key 2
\by B. C. Berndt\book Ramanujan's Notebooks
\bookinfo Part II
\publ Springer-Verlag
\publaddr New York \yr 1989
\endref

\ref\key 3
\by B. C. Berndt\book Ramanujan's Notebooks
\bookinfo Part III
\publ Springer-Verlag
\publaddr New York \yr 1991 
\endref

\ref\key 4
\by B. C. Berndt\paper Ramanujan's theory of theta-functions
\inbook Theta Functions From the Classical to the Modern
\ed M. Ram Murty
\bookinfo vol. 1 of CRM Proceedings \& Lecture Notes 
\publ American Mathematical Society\publaddr Providence, RI  
\yr 1993\pages 1--63
\endref

\ref\key 5
\by R. Borcherds
\moreref Private Communication, 9-14-2000
\endref

\ref\key 6
\by K. Chandrasekharan
\paper Elliptic Functions
\inbook vol. 281 of Grundlehren Math. Wiss.
\publ Springer-Verlag\publaddr Berlin
\yr 1985
\endref

\ref\key 7
\by L. Comtet\book Advanced Combinatorics
\publ D. Reidel Pub. Co.
\publaddr Dordrecht-Holland/Boston-U.S.A.\yr 1974
\endref

\ref\key 8
\by R. Dedekind
\paper Schreiben an Herrn Borchardt \"uber die Theorie der elliptischen
Modulfunktionen
\jour J. Reine Angew. Math. \vol 83\yr 1877 
\pages 265--292
\moreref\inbook reprinted in Richard Dedekind Gesammelte mathematische Werke, 
vol. 1, (Eds. R. Fricke, E. Noether, and \"O. Ore)
\publ Friedr.Vieweg \& Sohn Ake.-Ges.
\publaddr Braunschweig
\yr 1930\pages 174--201 
\moreref
\publ reprinted by Johnson Reprint Corporation
\publaddr New York, 1969
\endref

\ref\key 9
\by R. Dedekind
\paper Erl\"auterungen zu den Fragmenten XXVIII
\inbook Bernhard Riemann's Gesammelte mathematische Werke und wissenschaftlicher
Nachla{\ss}, Second Edition, (Eds. R. Dedekind and H. Weber)
\publ B. G. Teubner\publaddr Leipzig
\yr 1892\pages 466--478
\moreref 
\publ reprinted by Dover Publications, Inc.
\publaddr New York, 1953
\moreref\inbook also reprinted in Richard Dedekind Gesammelte mathematische Werke, 
vol. 1, (Eds. R. Fricke, E. Noether, and \"O. Ore)
\publ Friedr.Vieweg \& Sohn Ake.-Ges.
\publaddr Braunschweig
\yr 1930\pages 159--172 
\moreref
\publ reprinted by Johnson Reprint Corporation
\publaddr New York, 1969
\endref

\ref\key 10
\by R. Fricke
\paper Erl\"auterungen zur vorstehenden Abhandlung
\inbook Richard Dedekind Gesammelte mathematische Werke, 
vol. 1, (Eds. R. Fricke, E. Noether, and \"O. Ore)
\publ Friedr.Vieweg \& Sohn Ake.-Ges.
\publaddr Braunschweig
\yr 1930\page 173
\moreref
\publ reprinted by Johnson Reprint Corporation
\publaddr New York, 1969
\endref

\ref\key 11
\by F. Garvan
\moreref Private Communication, March 1997
\endref

\ref\key 12
\by I. S. Gradshteyn and I. M. Ryzhik
\book Table of Integrals, Series, and Products
\bookinfo 4th ed. \publ Academic Press
\publaddr San Diego \yr 1980 
\transl Translated from the Russian by Scripta Technica, Inc., and
edited by A. Jeffrey
\endref

\ref\key 13
\by H. Hancock
\book Lectures on the Theory of Elliptic Functions, Volume I Analysis
\publ John Wiley \& Sons
\publaddr New York\yr 1910 
\moreref 
\publ reprinted by Dover Publications, Inc.
\publaddr New York, 1958
\endref

\ref\key 14
\by A. Hurwitz
\paper Grundlagen einer independenten Theorie der elliptischen Modulfunctionen und
Theorie der Multiplicatorgleichungen erster Stufe
\paperinfo Ph. D. Dissertation
\publ University of Leipzig, 1881
\endref

\ref\key 15
\by A. Hurwitz
\paper Grundlagen einer independenten Theorie der elliptischen Modulfunctionen und
Theorie der Multiplicatorgleichungen erster Stufe
\jour Math. Ann.\vol 18\yr 1881 
\pages 528--592
\moreref 
reprinted in Mathematische Werke von Adolf Hurwitz;Herausgegeben von der Abteilung
f\"ur Mathematik und Physik der Eidgen\"ossischen Technischen Hochschule in
Z\"urich
\inbook\vol 1, Funktionentheorie
\publ E. Birkhauser \& cie.\publaddr Basel
\yr 1932\pages 1--66 
\endref  

\ref\key 16
\by C. G. J. Jacobi\paper Fundamenta Nova Theoriae
Functionum Ellipticarum\jour Regiomonti. Sumptibus
fratrum Borntr\"ager, 1829
\moreref reprinted in Jacobi's Gesammelte Werke
\inbook\vol 1\publ Reimer\publaddr Berlin 
\yr 1881--1891\pages 49--239 
\moreref 
\publ reprinted by Chelsea
\publaddr New York, 1969
\moreref 
Now distributed by The American Mathematical Society
\publaddr Providence, RI
\endref  

\ref\key 17
\by W. B. Jones and W. J. Thron\paper Continued Fractions:Analytic
Theory and Applications  
\inbook Encyclopedia of Mathematics and Its Applications\vol 11
\ed G.-C. Rota\publ Addison-Wesley\publaddr London
\yr 1980
\moreref
\publ Now distributed by Cambridge University Press  
\publaddr Cambridge
\endref

\ref\key 18
\by F. Klein and R. Fricke
\book Vorlesungen \"uber die Theorie der elliptischen Modulfunktionen
\bookinfo Vol. 1\publ B. G. Teubner\publaddr Leipzig
\yr 1890
\moreref
\publ reprinted by Johnson Reprint Corporation
\publaddr New York, 1966
\endref 

\ref\key 19
\by M. I. Knopp
\book Modular Functions in Analytic Number Theory
\publ Markham Publishing Company
\publaddr Chicago\yr 1970
\endref

\ref\key 20
\by C. Krattenthaler
\paper Advanced determinant calculus
\jour S\'eminaire Lotharingien de Combinatoire
\vol 42\yr 1999 \pages B42q (67 pp)
\endref

\ref\key 21
\by D. F. Lawden
\paper Elliptic Functions and Applications
\inbook vol. 80 of Applied Mathematical Sciences
\publ Springer-Verlag\publaddr New York 
\yr 1989
\endref

\ref\key 22
\by H. McKean and V. Moll
\book Elliptic Curves:function theory, geometry, arithmetic
\publ Cambridge University Press
\publaddr Cambridge\yr 1997
\endref

\ref\key 23
\by S. C. Milne\paper New infinite families of exact sums of
squares formulas, Jacobi elliptic functions, and Ramanujan's tau
function
\jour Proc. Nat. Acad. Sci., U.S.A.\vol 93\yr 1996
\pages 15004--15008
\endref

\ref\key 24
\by S. C. Milne\paper Infinite families of exact sums of squares formulas, 
Jacobi elliptic functions, continued fractions, and Schur functions
\paperinfo preprint arXiv:math.NT/0008068 (to appear in the Ramanujan Journal)
\endref

\ref\key 25
\by F. \`E. Molin (T. Molien)
\paper \"Uber die lineare Transformation der elliptischen Functionen
\paperinfo Master's Dissertation
\publ University of Derpt (Yuryev or Dorpat)-now Tartu, 1885
\endref

\ref\key 26
\by F. \`E. Molin (T. Molien)
\paper Ueber gewisse, in der Theorie der elliptischen Functionen auftretende
Einheitswurzeln
\jour Berichte \"uber die Verhandlungen der K\"onigl. S\"achsischen Geselischaft
der Wissenschaften zu Leipzig. Mathematisch-Physische Classe. Leipzig.
\vol 37\yr 1885\pages 25--38
\endref

\ref\key 27
\by S. Ramanujan\paper On certain arithmetical functions\jour Trans.
Cambridge Philos. Soc. \vol 22\yr 1916\pages 159--184
\moreref\inbook reprinted in Collected Papers of Srinivasa Ramanujan
\publ Chelsea\publaddr New York
\yr 1962\pages 136--162 
\moreref 
\publ reprinted by AMS Chelsea
\publaddr Providence, RI, 2000
\moreref 
Now distributed by The American Mathematical Society
\publaddr Providence, RI
\endref  

\ref\key 28
\by R. A. Rankin\book Modular Forms and Functions
\publ Cambridge University Press
\publaddr Cambridge \yr 1977 
\endref

\ref\key 29
\by B. Riemann
\paper Fragmente \"uber die Grenzf\"alle der elliptischen Modulfunctionen. (1852.)
\inbook Bernhard Riemann's Gesammelte mathematische Werke und wissenschaftlicher
Nachla{\ss}, Second Edition, (Eds. R. Dedekind and H. Weber)
\publ B. G. Teubner\publaddr Leipzig
\yr 1892\pages 455--465
\moreref 
\publ reprinted by Dover Publications, Inc.
\publaddr New York, 1953
\endref

\ref\key 30
\by B. Schoeneberg
\paper Elliptic Modular Functions
\inbook vol. 203 of Grundlehren Math. Wiss.
\publ Springer-Verlag\publaddr Berlin
\yr 1974
\endref

\ref\key 31
\by J.-P. Serre
\paper A Course in Arithmetic
\inbook vol. 7 of Graduate Texts in Mathematics
\publ Springer-Verlag\publaddr New York 
\yr 1973
\endref

\ref\key 32
\by N. J. A. Sloane and Simon Plouffe
\book The Encyclopedia of Integer Sequences
\publ Academic Press
\publaddr San Diego\yr 1995
\endref 

\ref\key 33
\by H. P. F. Swinnerton-Dyer
\paper Congruence properties of $\tau(n)$
\inbook Ramanujan Revisited\eds G. E. Andrews et al.
\publ Academic Press\publaddr New York\yr 1988 \pages 289--311\endref

\ref\key 34
\by K. Venkatachaliengar
\book Development of Elliptic Functions According to Ramanujan
\publ Tech. Report 2, Madurai Kamaraj University 
\publaddr Madurai \yr 1988
\endref

\ref\key 35
\by H. S. Wall
\book Analytic Theory of Continued Fractions
\publ D. Van Nostrand
\publaddr New York\yr 1948 
\moreref 
\publ reprinted by Chelsea
\publaddr New York, 1973
\endref

\ref\key 36
\by E. T. Whittaker and G. N. Watson
\book A Course of Modern Analysis\bookinfo 4th ed. 
\publ Cambridge University Press
\publaddr Cambridge \yr 1927 
\endref

\ref\key 37
\by S. Wolfram
\book The Mathematica Book
\bookinfo 4th ed. 
\publ Wolfram Media/Cambridge University Press
\publaddr Cambridge\yr 1999
\endref

\endRefs
\enddocument